\documentclass[12pt]{amsart}
\usepackage{amsmath}
\usepackage{amsxtra}
\usepackage{amscd}
\usepackage{amsthm}
\usepackage{amsfonts}
\usepackage{amssymb}
\usepackage{eucal}
\def\({\left(}
\def\){\right)}

\newcommand{\bea}{\begin{eqnarray}}
\newcommand{\ena}{\end{eqnarray}}
\newcommand{\be}{\begin{eqnarray*}}
\newcommand{\en}{\end{eqnarray*}}
\newcommand{\ba}{\begin{array}}
\newcommand{\ea}{\end{array}}

\newcommand{\Z}{{\mathbb Z}}

\newenvironment{tenumerate}{
  \begin{enumerate}
  
  }{\end{enumerate}}
\newcommand{\bi}{\begin{tenumerate}}
\newcommand{\ei}{\end{tenumerate}}
\newcommand{\isoto}[1][]%
{{\mathop{\buildrel{\sim}\over\longrightarrow}\limits_{#1}}}

\newcommand{\AHA}{\mathcal{H}^{\mathrm{aff}}_n}

\numberwithin{equation}{section}
\newtheorem{thm}{Theorem}[section]
\newtheorem{prop}[thm]{Proposition}
\newtheorem{lem}[thm]{Lemma}
\newtheorem{rem}[thm]{Remark}

\newtheorem{definition}[thm]{Definition}

\newtheorem{defn}[thm]{Definition}

\begin{document}

\title[The qKZ equation and non-symmetric Macdonald polynomials]
{The quantum Knizhnik-Zamolodchikov equation and 
non-symmetric Macdonald polynomials}
\author{M.~Kasatani and Y.~Takeyama}
\address{MK: Department of Mathematics, Graduate School of Science,
Kyoto University, Kyoto 606-8502, 
Japan}\email{kasatani@math.kyoto-u.ac.jp}
\address{YT: Department of Mathematics, 
Graduate School of Pure and Applied Sciences, 
Tsukuba University, Tsukuba, Ibaraki 305-8571, Japan}
\email{takeyama@math.tsukuba.ac.jp}

\begin{abstract}
We construct special solutions of the quantum Knizhnik-Zamolodchikov equation 
on the tensor product of the vector representation of 
the quantum algebra of type $A_{N-1}$. 
They are constructed from non-symmetric Macdonald polynomials 
through the action of the affine Hecke algebra. 
As a special case the matrix element of the vertex operators of level one is reproduced. 
\end{abstract}

\maketitle

{\small
\noindent{\it 2000 Math. Subj. Class.} 
39A13, 33C52, 81R50.
}
\medskip 

\setcounter{section}{0}
\setcounter{equation}{0}

\section{Introduction}

In this paper we construct 
special solutions of the quantum Knizhnik-Zamolodchikov (qKZ) equation 
in terms of non-symmetric Macdonald polynomials.  

The qKZ equation, derived by Frenkel and Reshetikhin \cite{FR}, 
is the system of difference equations satisfied by 
matrix elements of the vertex operators 
in the representation theory of the quantum affine algebra. 
In this paper we consider the qKZ equation 
on the tensor product of the vector representation of 
the quantum algebra $U_{q}(sl_{N})$: 
\begin{eqnarray*}
&& 
G(z_{1}, \ldots , pz_{m}, \ldots , z_{n})=
R_{m, m-1}(pz_{m}/z_{m-1}) \cdots R_{m, 1}(pz_{m}/z_{1})\, 
{\textstyle ( \prod_{j=1}^{N-1}\kappa_{j}^{h_{j}} )}_{m} \\ 
&& \qquad  {}\times 
R_{n,m}(z_{n}/z_{m})^{-1} \cdots R_{m+1, m}(z_{m+1}/z_{m})^{-1} \,
G(z_{1}, \ldots , z_{m}, \ldots , z_{n}).  
\end{eqnarray*} 
Here $G(z_{1}, \ldots , z_{n})$ is an unknown function taking values in $V^{\otimes n}$, 
where $V \simeq \mathbb{C}^{N}$ is the vector representation. 
The operator $R(z)$ is the $R$-matrix (see \eqref{eq:def-Rmatrix} below), 
$h_{j} \, (j=1, \ldots , N-1)$ is the basis of the Cartan subalgebra of $sl_{N}$, 
and $p, \kappa_{1}, \ldots , \kappa_{N-1}$ are parameters of the equation. 
The indices of the operators in the right hand side signify 
the positions of the components in $V^{\otimes n}$ where the operator acts. 
The value $\ell$ determined by $p=q^{2(N+\ell)}$ is called a level of the qKZ equation.

There are some methods to construct solutions of the qKZ equation. 
One of them is to use multiple integrals of the hypergeometric type \cite{Mi, VT, MTT}. 
This works for any parameters $p, \kappa_{1}, \ldots , \kappa_{N-1}$ 
such that the multiple integrals converge. 
Another method is the bosonization of vertex operators. 
For the integrable irreducible highest weight $U_{q}(\widehat{sl}_{N})$-modules of level one, 
the bosonization is constructed by Koyama \cite{Ko}. 
By using it, Nakayashiki calculated the matrix element of the vertex operators \cite{N}. 
It gives by definition a solution to the qKZ equation of level one, 
where the parameters $\kappa_{1}, \ldots , \kappa_{N-1}$ are determined from 
the highest weight. 

Recently Di Francesco and Zinn-Justin constructed a polynomial solution 
in the case of level one \cite{DZ} by using the representation theory 
of the affine Hecke algebra (AHA).
In a similar manner Kasatani and Pasquier obtained a solution of the qKZ equation 
of level $-1/2$ associated with $U_{q}(\widehat{sl}_{2})$ \cite{KP}. 
In this paper we generalize these results to the case of $U_{q}(\widehat{sl}_{N})$ and other levels. 


Let us give a sketch of our construction of solutions. 
We use the spin basis instead of the path basis in the construction of \cite{DZ} and \cite{KP}. 
Expand the unknown function $G(z_{1}, \ldots , z_{n})$ into a linear combination of 
the tensor products $v_{\epsilon_{1}} \otimes \cdots \otimes v_{\epsilon_{n}}$, 
where $\{v_{\epsilon}\}_{\epsilon=0}^{N-1}$ is the standard basis of $V$. 
We consider the set of functions which appear in the expansion as coefficients. 
The qKZ equation can be described as a condition of constraint for the functions. 
In this paper we consider a stronger condition than the qKZ equation itself, 
and call a set of functions satisfying the condition a {\it qKZ family} 
(see Definition \ref{def:qKZ-family} below). 

The defining condition of qKZ family is described in terms of the action of the AHA 
on the space of functions. 
The generators of the AHA consist of two sets of elements 
$T_{i} \, (1 \le i \le n-1)$ and $Y_{j} \, (1 \le j \le n)$ satisfying 
some relations (see Definition \ref{def:AHA} below). 
The actions of $T_{i}$ and $Y_{j}$ are given by 
the Demazure-Lusztig operator and the $q$-Dunkl operator, respectively. 
{}From a viewpoint of the representation theory, 
a qKZ family is a set of vectors which move to each other 
by the action of the generators of the AHA. 
Moreover, if one vector is known, then all the vectors are determined 
through the action of the AHA. 
Hence the linear span of the vectors of a qKZ family 
determines a cyclic module of the AHA. 

Now we return to the description of our construction of solutions. 
{}From the definition, a qKZ family contains 
a joint eigenfunction of the $q$-Dunkl operators $Y_{j}$. 
Moreover, it is also an eigenfunction of some of 
the Demazure-Lusztig operators $T_{i}$.  
We prove that such an eigenfunction, conversely, generates a qKZ family. 
Thus construction of a qKZ family is reduced to that of 
an eigenfunction of the $q$-Dunkl operators and the Demazure-Lusztig operators. 
As is well known, 
non-symmetric Macdonald polynomials \cite{Ma} are such eigenfunctions.  
Therefore we can construct special solutions of the qKZ equation from 
non-symmetric Macdonald polynomials. 

Cherednik \cite{C} and Kato \cite{Kat} unveiled 
the relation between the qKZ equation and 
the eigenvalue problem of the Macdonald type:   
a certain linear combination of the coefficients in a solution of the qKZ equation 
gives a symmetric joint eigenfunction of the $q$-Dunkl operators 
(see \cite{Mi} for the explicit formula in the case of $n=N$). 
Our construction is consistent with this result because 
symmetric Macdonald polynomials can be obtained as a linear combination of non-symmetric ones. 

The plan of this paper is as follows. 
First we recall the definition of the qKZ equation in Section \ref{sec:qKZ-def}. 
In Section \ref{sec:qKZ-family} we give the definition of qKZ family, 
and prove that a qKZ family is constructed from 
a joint eigenfunction of the $q$-Dunkl operators and some of the Demazure-Lusztig operators. 
In Section \ref{sec:non-symm-Mac} we give explicitly the construction explained above 
of solutions of the qKZ equation 
in the case where the level is generic or a value of the form 
$\frac{k+1}{r-1}-N$, 
where $k$ and $r$ are positive integers such that 
$1\leq k\leq {\rm min}\{n-1, N\}, \, r \ge 2$, and $k+1$ and $r-1$ are coprime.
Then the parameters $\kappa_{1}, \ldots , \kappa_{N-1}$ are determined from 
the eigenvalues of the non-symmetric Macdonald polynomial 
for the $q$-Dunkl operators. 
Here it should be noted that in the latter case we need to specialize 
the two parameters in non-symmetric Macdonald polynomials, 
some of which are proved to be well-defined in \cite{Kas}. 
In Section \ref{sec:VO} we see that 
the matrix element of the vertex operators of level one is reproduced 
by the construction in Section \ref{sec:non-symm-Mac} 
in the case where $k=N$ and $r=2$.

\section{The quantum Knizhnik-Zamolodchikov equation}\label{sec:qKZ-def}

Let $V=\oplus_{\epsilon=0}^{N-1}\mathbb{C}\,v_{\epsilon}$ be the $N$-dimensional vector space. 
We regard $V$ as the vector representation of the quantum algebra $U_{q}(sl_{N})$. 
Define the linear operator $\bar{R}(z)$ acting on $V^{\otimes 2}$ by 
\begin{eqnarray*}
\bar{R}(z) \, (v_{\epsilon_{1}} \otimes v_{\epsilon_{2}})=
\sum_{\epsilon_{1}', \epsilon_{2}'}\bar{R}(z)_{\epsilon_{1}'\epsilon_{2}'}^{\epsilon_{1}\epsilon_{2}}\,
v_{\epsilon_{1}'}\otimes v_{\epsilon_{2}'},  
\end{eqnarray*}
where 
\begin{eqnarray*}
\bar{R}(z)_{ii}^{ii}=1, \quad 
\bar{R}(z)_{ij}^{ij}=\frac{(1-z)q}{1-q^2 z}, \quad 
\bar{R}(z)_{ij}^{ji}=\frac{1-q^2}{1-q^2 z}z^{\theta(i>j)} \quad (i\not=j)   
\end{eqnarray*}
and $\bar{R}(z)_{i'j'}^{ij}=0$ otherwise. 
Here 
\begin{eqnarray}
\theta(P)=\left\{
\begin{array}{ll}
1 & \hbox{if $P$ is true}, \\
0 & \hbox{if $P$ is false}.
\end{array}
\right.
\label{eq:def-theta}
\end{eqnarray}
Throughout this paper we assume that 
\begin{eqnarray*}
0<q<1.  
\end{eqnarray*}
Then the $R$-matrix $R(z)$ is given as follows \cite{DO}:
\begin{eqnarray}
R(z):=r(z)\bar{R}(z). 
\label{eq:def-Rmatrix} 
\end{eqnarray}
Here $r(z)$ is the normalization factor 
\begin{eqnarray*}
r(z)=q^{\frac{1}{N}-1}
\frac{(q^2 z ; q^{2N})_{\infty} (q^{2N-2} z ; q^{2N})_{\infty}}
{(z ; q^{2N})_{\infty} (q^{2N} z ; q^{2N})_{\infty}}, 
\qquad 
(z; x)_{\infty}:=\prod_{j=0}^{\infty}(1-x^{j}z). 
\end{eqnarray*} 
The matrix $R(z)$ is nothing but the image in ${\rm End}(V^{\otimes 2})$ 
of the universal $R$-matrix $\mathcal{R}'(z)$ of 
the quantum affine algebra $U_{q}(\widehat{sl}_{N})$ 
in the sense of Appendix 1 in \cite{IIJMNT}.

The qKZ equation is the following system of difference equations 
for an unknown function $G(z_{1}, \ldots , z_{n})$ taking values in $V^{\otimes n}$: 
\begin{eqnarray}
&& 
G(z_{1}, \ldots , pz_{m}, \ldots , z_{n})=
R_{m, m-1}(pz_{m}/z_{m-1}) \cdots R_{m, 1}(pz_{m}/z_{1})\, 
{\textstyle ( \prod_{j=1}^{N-1}\kappa_{j}^{h_{j}} )_{m}} \\ 
&& \hspace{5em} {}\times 
R_{n,m}(z_{n}/z_{m})^{-1} \cdots R_{m+1, m}(z_{m+1}/z_{m})^{-1} \,
G(z_{1}, \ldots , z_{m}, \ldots , z_{n}) 
\nonumber  
\end{eqnarray}
for $m=1, \ldots , n$. 
Here $R_{m, l}(z)$ is the operator acting on the tensor product 
of the $m$-th and the $l$-th components in $V^{\otimes n}$ as the $R$-matrix $R(z)$
\footnote{Note the order of indices: 
$R_{21}(z)=P\,R(z)\,P \not=R_{12}(z)$, where $P$ is the transposition $P(u \otimes v):=v \otimes u$}.  
The operator $(\prod_{j}\kappa_{j}^{h_{j}})_{m}$ acts on the $m$-th component, 
where $h_{j} \, (j=1, \ldots , N-1)$ is the basis of the Cartan subalgebra of $sl_{N}$.  
The action of $h_{j}$ on $V$ is given by 
\begin{eqnarray*}
h_{j}v_{j-1}=v_{j-1}, \quad h_{j}v_{j}=-v_{j}, \quad 
h_{j}v_{i}=0 \,\, (i\not=j-1, j). 
\end{eqnarray*} 
The complex numbers $\kappa_{1}, \ldots , \kappa_{N-1}$ are parameters of the qKZ equation. 
For the sake of simplicity, hereafter we assume that 
the difference step $p$ is a positive real number. 
When $p=q^{2(N+\ell)}$ the number $\ell$ is called {\it level}.

\section{qKZ family}\label{sec:qKZ-family}

\subsection{Affine Hecke algebra}

Let us summarize the basic facts about the affine Hecke algebra. 
We use the notation in \cite{MN}. 

\begin{definition}\label{def:AHA}
The affine Hecke algebra $\AHA$ of type $GL_{n}$ is an associative 
$\mathbb{C}(t^{1/2})$-algebra generated by $T_{i} \, (i=1, \ldots , n-1)$ and 
$Y_{j} \, (j=1, \ldots , n)$ satisfying the following relations: 
\begin{eqnarray*}
&& 
(T_{i}-t^{1/2})(T_{i}+t^{-1/2})=0 \quad (1 \le i \le n-1), \\ 
&& 
T_{i}T_{i+1}T_{i}=T_{i+1}T_{i}T_{i+1} \quad (1 \le i \le n-2), \\ 
&& 
T_{i}T_{j}=T_{j}T_{i} \quad (|i-j|>1), \\ 
&& 
Y_{i}Y_{j}=Y_{j}Y_{i} \quad (1 \le i, j \le n), \\ 
&& 
Y_{i}T_{j}=T_{j}Y_{i} \quad (j \not=i-1, i), \\ 
&& 
T_{i}Y_{i+1}T_{i}=Y_{i} \quad (1 \le i \le n-1). 
\end{eqnarray*}
\end{definition}

Define $\sigma \in \AHA$ by 
\begin{eqnarray*}
\sigma:= T_{n-1}^{-1} \cdots T_{i}^{-1}Y_{i} T_{i-1} \cdots T_{1}. 
\end{eqnarray*}
Note that the right hand side above does not depend on the value $i$. 
Then it is easy to see that $\sigma^{n}$ is central and 
$\sigma T_{i}=T_{i-1} \sigma \, (1<i<n)$. 
Moreover the algebra $\AHA$ is generated by 
$T_{i} \, (i=1, \ldots , n-1)$ and $\sigma$. 

Denote the Laurent polynomial ring with $n$ variables by 
$P_{n}=\mathbb{C}[z_{1}^{\pm 1}, \ldots , z_{n}^{\pm 1}]$. 
Let $\widehat{T}_{i} \, (i=1, \ldots , n-1)$ and 
$\omega$ be the linear operators on $P_{n}$ defined by 
\begin{eqnarray}
&& 
\widehat{T}_{i}:=t^{1/2}\tau_{i}+\frac{t^{1/2}-t^{-1/2}}{z_{i}/z_{i+1}-1}(\tau_{i}-1), 
\label{eq:action-1} \\ 
&& 
(\omega f)(z_{1}, \ldots , z_{n}):=f(pz_{n}, z_{1}, \ldots , z_{n-1}).
\label{eq:action-2}
\end{eqnarray}
Here $\tau_{i}$ is the permutation of the variables $z_{i}$ and $z_{i+1}$, 
and $p$ is a parameter. 
The operator $\widehat{T}_{i}$ is called the {\it Demazure-Lusztig operator}. 
We will identify the parameter $p$ with the difference step $p$ 
in the qKZ equation. 

\begin{prop}
The linear map $\pi: \AHA \to {\rm End}(P_{n})$ defined by 
$\pi(T_{i})=\widehat{T}_{i} \, (i=1, \ldots , n-1)$ and $\pi(\sigma)=\omega$ 
gives a representation of $\AHA$. 
\end{prop}

\subsection{qKZ family}\label{subsec:qKZ-family}

Hereafter we assume that 
\begin{eqnarray*}
n \ge N \ge 2. 
\end{eqnarray*}
Let $d_{0}, \ldots , d_{N-1}$ be positive integers satisfying 
$\sum_{j=0}^{N-1}d_{j}=n$. 
Denote by $I_{d_{0}, \ldots , d_{N-1}}$ the set of $n$-tuples 
$\epsilon=(\epsilon_{1}, \ldots , \epsilon_{n})$ satisfying 
\begin{eqnarray*}
\#\{a \, | \, \epsilon_{a}=j \}=d_{j} \qquad (0 \le j \le N-1).
\end{eqnarray*}

Now we give the definition of qKZ family: 
\begin{definition}\label{def:qKZ-family}
A set of Laurent polynomials 
\begin{eqnarray*}
\{ f_{\epsilon_{1}, \ldots , \epsilon_{n}} \in P_{n} \, | \, 
 (\epsilon_{1}, \ldots , \epsilon_{n}) \in I_{d_{0}, \ldots , d_{N-1}} \} 
\end{eqnarray*} 
is called a qKZ family 
of sign $(\pm)$ with exponents $(c_{0}, \ldots , c_{N-1})$ 
if it satisfies the following conditions:
\begin{itemize}
\item If $\epsilon_{i}=\epsilon_{i+1}$, then
$\widehat{T}_{i}f_{\ldots , \epsilon_{i}, \epsilon_{i+1}, \ldots }=
\pm t^{\pm 1/2} f_{\ldots , \epsilon_{i}, \epsilon_{i+1}, \ldots }.$

\item If $\epsilon_{i}>\epsilon_{i+1}$, then
$\widehat{T}_{i}f_{\ldots , \epsilon_{i}, \epsilon_{i+1}, \ldots }=
f_{\ldots , \epsilon_{i+1}, \epsilon_{i}, \ldots }.$

\item 
$\omega f_{\epsilon_{n}, \epsilon_{1}, \ldots , \epsilon_{n-1}}=
c_{\epsilon_{n}}f_{\epsilon_{1}, \ldots , \epsilon_{n}}.$
\end{itemize}
Here the operators $\widehat{T}_{i} \, (i=1, \ldots , n-1)$ and $\omega$ are 
defined by \eqref{eq:action-1} and \eqref{eq:action-2}, respectively. 
\end{definition}

Note that a qKZ family $\{f_{\epsilon_{1}, \ldots , \epsilon_{n}}\}$
is uniquely determined from one member of it through the action of $\widehat{T}_{i}$. 
Hence the linear span $\sum \mathbb{C}f_{\epsilon_{1}, \ldots , \epsilon_{n}}$ is a 
cyclic $\AHA$-submodule of $P_{n}$.

In the rest of this subsection 
we show that a solution of the qKZ equation can be constructed 
{}from a qKZ family. 
 
Let ${\bf f}=\{f_{\epsilon_{1}, \ldots , \epsilon_{n}}\}$ be a qKZ family 
with exponents $(c_{0}, \ldots , c_{N-1})$. 
Now we determine two parameters $\alpha$ and $\beta$, and a function $h(z)$ 
according to the sign of ${\bf f}$ as follows. 
If the sign of ${\bf f}$ is plus, 
we define $\alpha, \beta$ by 
\begin{eqnarray*}
p^{\alpha}=(\prod_{j=0}^{N-1}c_{j})^{-1/N} \, q^{-(n+1)(1/N-1)}, \qquad 
p^{\beta}=q^{2(1/N-1)}, 
\end{eqnarray*}
and take a solution $h(z)$ of the difference equation 
\begin{eqnarray*}
\frac{h(p^{-1}z)}{h(z)}=
\frac{(z ; q^{2N})_{\infty} (q^{2N}z; q^{2N})_{\infty}}
{(q^{2}z; q^{2N})_{\infty} (q^{2N-2}z; q^{2N})_{\infty}}.
\end{eqnarray*}
Similarly, in the case where the sign of ${\bf f}$ is minus, 
we determine $\alpha, \beta$ and $h(z)$ by the following formulas: 
\begin{eqnarray}
&& 
p^{\alpha}=(-1)^{n-1}\,(\prod_{j=0}^{N-1}c_{j})^{-1/N} \, q^{-(n+1)(1+1/N)}, \qquad 
p^{\beta}=q^{2(1+1/N)}, 
\label{eq:def-alphabeta}\\
&&
\frac{h(p^{-1}z)}{h(z)}=
\frac{(z ; q^{2N})_{\infty} (q^{2N}z; q^{2N})_{\infty}}
{(q^{2(N+1)}z; q^{2N})_{\infty} (q^{-2}z; q^{2N})_{\infty}}.
\label{eq:def-h}
\end{eqnarray}

Now let us construct a solution of the qKZ equation. 
Define the function $K(z_{1}, \ldots , z_{n})$ by 
\begin{eqnarray}
K(z_{1}, \ldots , z_{n}):=\prod_{a=1}^{n}z_{a}^{\alpha+\beta a} 
\prod_{1 \le a<b \le n}h(z_{b}/z_{a}). 
\label{eq:normalization-factor}
\end{eqnarray}
and the $V^{\otimes n}$-valued function $F(z_{1}, \ldots , z_{n})$ by 
\begin{eqnarray*}
F(z_{1}, \ldots , z_{n}):=
\sum_{(\epsilon_{1}, \ldots , \epsilon_{n}) \in I_{d_{0}, \ldots , d_{N-1}}} 
f_{\epsilon_{1}, \ldots , \epsilon_{n}}(z_{1}, \ldots , z_{n}) \, 
v_{\epsilon_{1}} \otimes \cdots \otimes v_{\epsilon_{n}}.
\end{eqnarray*}
Set 
\begin{eqnarray}
G(z_{1}, \ldots , z_{n}):=K(z_{1}, \ldots , z_{n})F(z_{1}, \ldots , z_{n}).
\label{eq:qKZ-solution}
\end{eqnarray}

\begin{prop}\label{prop:qKZ-solution}
Let ${\bf f}=\{f_{\epsilon_{1}, \ldots , \epsilon_{n}}\}$ be a qKZ family 
of sign $(\pm)$ with exponents $(c_{0}, \ldots , c_{N-1})$. 
Then $G(z_{1}, \ldots , z_{n})$ is a solution of the qKZ equation whose 
parameters $q$ and $\kappa_{j} \, (j=1, \ldots , N-1)$ are determined by $q=\pm t^{\pm 1/2}$ and 
\begin{eqnarray}
\kappa_{j}=\prod_{l=0}^{j-1}c_{l} \cdot (\prod_{l=0}^{N-1}c_{l})^{-j/N}.
\label{eq:def-kappa}
\end{eqnarray} 
\end{prop}

\begin{proof}
Here we give the proof in the case where the sign of ${\bf f}$ is minus. 
The proof for the case of plus sign is similar. 

Suppose that two functions $F$ and $G$ are related by \eqref{eq:qKZ-solution}. 
Then the qKZ equation for $G$ is equivalent to 
the following equation for $F$:
\begin{eqnarray}
&& 
F(z_{1}, \ldots , pz_{m}, \ldots , z_{n}) 
\label{eq:qKZ-polynomial-form} \\ 
&& {}=(-1)^{n-1}
q^{-2(n-2m+1)} \, (\prod_{j=0}^{N-1}c_{j})^{1/N}
\prod_{a=1}^{m-1}
\frac{1-q^{-2}z_{a}/pz_{m}}{1-q^{2}z_{a}/pz_{m}} \, 
\prod_{b=m+1}^{n}
\frac{1-q^{2}z_{m}/z_{b}}{1-q^{-2}z_{m}/z_{b}} 
\nonumber \\ 
&& {}\times 
\bar{R}_{m, m-1}(pz_{m}/z_{m-1}) \cdots \bar{R}_{m, 1}(pz_{m}/z_{1})\, 
( \prod_{j=1}^{N-1}\kappa_{j}^{h_{j}} )_{m} 
\nonumber \\ 
&& {}\times 
\bar{R}_{n,m}(z_{n}/z_{m})^{-1} \cdots \bar{R}_{m+1, m}(z_{m+1}/z_{m})^{-1} \,
F(z_{1}, \ldots , z_{m}, \ldots , z_{n}). 
\nonumber 
\end{eqnarray}

On the other hand, 
if ${\bf f}$ is a qKZ family of sign $(-)$ with exponents $(c_{0}, \ldots , c_{N-1})$, 
we obtain the following equalities by setting $t^{1/2}=-q^{-1}$:
\begin{eqnarray*}
&& 
F( \ldots , z_{i+1}, z_{i}, \ldots ) \\ 
&& {}=
(-q^{-2}) \frac{1-q^2 z_{i}/z_{i+1}}{1-q^{-2}z_{i}/z_{i+1}}
P_{i, i+1}\bar{R}_{i, i+i}(z_{i}/z_{i+1}) 
F(\ldots , z_{i}, z_{i+1}, \ldots ), 
\end{eqnarray*}
where $P$ is the transposition $P(u \otimes v):=v \otimes u$, and 
\begin{eqnarray*}
P_{n-1,n}\cdots P_{1,2}F(pz_{n}, z_{1}, \ldots , z_{n-1})=
(\prod_{j=0}^{N-1}c_{j})^{1/N} \, 
(\prod_{j=1}^{N-1}\kappa_{j}^{h_{j}})_{n}\, 
F(z_{1}, \ldots , z_{n}).
\end{eqnarray*}
Here parameters $\kappa_{j} \, (j=1, \ldots , N-1)$ are defined by \eqref{eq:def-kappa}. 
It is easy to derive \eqref{eq:qKZ-polynomial-form} from the equalities above 
and $\bar{R}_{12}(z)^{-1}=\bar{R}_{21}(z^{-1})$. 
\end{proof}

\subsection{Equivalence to the eigenvalue problem}

Hereafter we often use the short notation $\epsilon=(\epsilon_{1}, \ldots , \epsilon_{n})$ 
to specify an element of $I_{d_{0}, \ldots , d_{N-1}}$. 
Let $\{f_{\epsilon}\}$ be a qKZ family of sign $(\pm)$. 
Consider the member $f_{\delta}$, where 
\begin{eqnarray*}
\delta:=(0^{d_0},1^{d_1},\cdots, (N-1)^{d_{N-1}}). 
\end{eqnarray*}
Then it satisfies $\widehat{T}_{i}f_{\delta}=\pm t^{\pm1/2}f_{\delta}$ 
for $1 \le i \le n-1$ such that $\delta_{i}=\delta_{i+1}$. 
Moreover it is an eigenfunction of the {\it $q$-Dunkl operators} 
\begin{eqnarray*}
\widehat{Y}_{j}:=\pi(Y_{j})=
\widehat{T}_{j} \cdots \widehat{T}_{n-1} \omega 
\widehat{T}_{1}^{\, -1} \cdots \widehat{T}_{j-1}^{\, -1}.
\end{eqnarray*}
Thus a qKZ family contains a joint eigenfunction 
of the $q$-Dunkl operators $\widehat{Y}_{j}$ and 
some of the Demazure-Lusztig operators $\widehat{T}_{i}$. 
On the contrary we can construct a qKZ family from such an eigenfunction as follows. 

Let us introduce some notation.  
An element $\lambda=(\lambda_{1}, \ldots , \lambda_{n}) \in \mathbb{Z}^{n}$ is called 
{\it dominant} (or {\it anti-dominant}) if $\lambda_{1} \ge \cdots \ge \lambda_{n}$ 
(or $\lambda_{1} \le \cdots \le \lambda_{n}$, resp.). 
The symmetric group $S_{n}$ acts on $\mathbb{Z}^{n}$ by 
$\sigma \lambda:=(\lambda_{\sigma^{-1}(1)}, \ldots , \lambda_{\sigma^{-1}(n)})$. 
We denote the orbit of $\lambda \in \mathbb{Z}^{n}$ by $S_{n}\lambda$.  

\begin{definition}\label{def:shortest-element}
For $\lambda\in\Z^n$, we denote by 
$\lambda^+$ $(\lambda^-)$ the unique dominant {\rm (}anti-dominant{\rm )} element in $S_n\lambda$,
respectively.
We denote by $w_\lambda^+$ $(w_\lambda^-)$ the shortest element in $S_n$ such
that $w_\lambda^+\lambda^+=\lambda$ $(w_\lambda^-\lambda=\lambda^-)$,
respectively.
\end{definition}

Since $I_{d_{0}, \ldots, d_{N-1}}$ is a subset of $\mathbb{Z}^{n}$,  
we use the notation above also for the elements of $I_{d_{0}, \ldots , d_{N-1}}$. 
For example, we have $w_{\epsilon}^{-}\epsilon=\delta$ 
for any $\epsilon \in I_{d_{0}, \ldots , d_{N-1}}$. 

For $w \in S_{n}$ we denote its length by $\ell(w)$.  
Let $w=s_{i_{1}} \cdots s_{i_{m}}$ be a reduced expression, 
where $s_{i}$ is the transposition $s_{i}=(i, i+1)$. 
Then we set $\widehat{T}_{w}:=\widehat{T}_{i_{1}} \cdots \widehat{T}_{i_{m}}$. 
This does not depend on the choice of reduced expression of $w$.  

Now we are in position to prove the main theorem 
which plays a key role in the next section: 

\begin{thm}\label{thm:constr}
Fix positive integers $d_0, \ldots, d_{N-1}$ satisfying 
$\sum_{j=0}^{N-1}d_{j}=n$ 
and set $\delta=(0^{d_0}, 1^{d_1}, \cdots, (N-1)^{d_{N-1}})$.
Suppose that $E=E(z_{1}, \ldots , z_{n})$ is a solution to 
the following eigenvalue problem:
\begin{eqnarray}
\widehat{Y}_j E&=&\chi_{j}E \quad (1 \leq \forall{j} \leq n)
\label{eq:equiv1}\\
\widehat{T}_i E&=& \pm t^{\pm 1/2}E \quad \hbox{if \, $\delta_i=\delta_{i+1}$}
\label{eq:equiv2}.
\end{eqnarray}
Here the sign in the right hand side \eqref{eq:equiv2} should be independent on $i$.
Set $f_{\epsilon}:=(\widehat{T}_{w_{\epsilon}^{-}})^{-1}E$ for 
$\epsilon \in I_{d_{0}, \ldots , d_{N-1}}$.  
Then $\{ f_{\epsilon}\}$ is a qKZ family of sign $(\pm)$ with exponents 
$c_i=\chi_{d_0+\cdots+d_i}(\pm t^{\pm1/2})^{d_i-1} \, (0 \le i \le N-1)$.
\end{thm}

\begin{rem}
The consistency of the eigenvalue problem \eqref{eq:equiv1} and \eqref{eq:equiv2} 
implies that 
the eigenvalues $\chi_{i}$ should satisfy 
$\chi_i=t^{\pm1}\chi_{i+1}$ if $\delta_i=\delta_{i+1}$. 
Hence all the eigenvalues $\chi_i$ are restored from the exponents $c_{i}$.  
\end{rem}

\begin{proof}[Proof of Theorem \ref{thm:constr}]

Let $E$ be a solution to the eigenvalue problem 
\eqref{eq:equiv1} and \eqref{eq:equiv2}, 
and set $f_{\epsilon}=(\widehat{T}_{w_{\epsilon}^{-}})^{-1}E$. 
Note that $f_{\delta}=E$. 
Let us check that the family of the functions 
$\{ f_{\epsilon} \}_{\epsilon \in I_{d_0,\ldots,d_{N-1}}}$
satisfies the three properties in Definition \ref{def:qKZ-family} in order. 

Take $\epsilon\in I_{d_0, \ldots, d_{N-1}}$.
Suppose $\epsilon_i=\epsilon_{i+1}$ and let $j=w_\epsilon^-(i)$.
Since $w_\epsilon^-$ is taken to be shortest, we see that
$\delta_j=\delta_{j+1}$,
$w_\epsilon^- s_i=s_j w_\epsilon^-$ and
$\ell(w_\epsilon^- s_i)=\ell(s_j w_\epsilon^-)=\ell(w_\epsilon^-)+1$.
Hence
\begin{eqnarray*}
\widehat{T}_{w_\epsilon^-}\widehat{T}_if_\epsilon
&=&\widehat{T}_j\widehat{T}_{w_\epsilon^-}f_\epsilon \\
&=&\widehat{T}_jf_\delta \\
&=&(\pm t^{\pm1/2})f_\delta.
\end{eqnarray*}
By applying $(\widehat{T}_{w_\epsilon^-})^{-1}$ on both sides, we obtain
the first property.

If $\epsilon_i>\epsilon_{i+1}$, then
$w_{\epsilon}^-=w_{s_i\epsilon}^-s_i$ and $\ell(w_\epsilon^-)=\ell(w_{s_i\epsilon}^-)+1$.
Hence
\begin{eqnarray*}
\widehat{T}_{w_{s_i\epsilon}^-}f_{s_i\epsilon}&=&f_\delta\\
&=&\widehat{T}_{w_{\epsilon}^-}f_\epsilon \\
&=&\widehat{T}_{w_{s_i\epsilon}^-}\widehat{T}_if_{\epsilon}.
\end{eqnarray*}
This gives the second property.

Let us check the third property. 
Let $\epsilon\in I_{d_0, \ldots, d_{N-1}}$ and set 
\begin{eqnarray*}
i_0:=w_\epsilon^-(n)=d_{0}+\cdots +d_{\epsilon_{n}}. 
\end{eqnarray*}
Then $\delta_{i_0}=\epsilon_n$.
Since $\sharp\{1\leq i<i_0 \, | \,  \delta_i=\delta_{i_0}\}=d_{\epsilon_{n}}-1$, we have 
\begin{eqnarray}
f_\delta&=&\chi_{i_0}^{-1}\widehat{T}_{i_0}\cdots \widehat{T}_{n-1}\omega
\widehat{T}_1^{-1}\cdots \widehat{T}_{i_0-1}^{-1}f_\delta \label{eq:constr1} \\
&=&(\pm t^{\pm1/2})^{-d_{\epsilon_{n}}+1}\chi_{i_0}^{-1}\widehat{T}_{i_0}\cdots \widehat{T}_{n-1}
\omega 
f_{\delta_{i_0},\delta_1,\ldots,\hat{\delta_{i_0}},\ldots,\delta_n}.
\nonumber
\end{eqnarray}
Set $\delta'=(\delta_{i_0},\delta_1,\ldots,\hat{\delta_{i_0}},\ldots,\delta_n)$
and $\epsilon'=(\epsilon_n,\epsilon_1,\ldots,\epsilon_{n-1})$.
Let $w_1$ be the shortest element such that
$w_1\epsilon'=\delta'$.
Then
$w_{\epsilon'}^-=w_{\delta'}^-w_1$
 and $\ell(w_{\epsilon'}^-)=\ell(w_{\delta'}^-)+\ell(w_1)$.
Hence 
\bea
f_{\delta'}=\widehat{T}_{w_1}f_{\epsilon'}. \label{eq:constr2}
\ena
On the other hand,
let $w_2$ be the shortest element such that
$w_2\epsilon=(\delta_1,\ldots,\hat{\delta_{i_0}},\ldots,\delta_n,\delta_{i_0})$.
Since there is no $i$ such that $i_0< i\leq n$ and $\delta_i=\delta_{i_0}$,
we see that
$w_\epsilon^-=s_{i_0}\cdots s_{n-1}w_2$ and $\ell(w_\epsilon^-)=n-i_0+\ell(w_2)$.
Therefore 
\bea
f_\delta&=&\widehat{T}_{w_\epsilon^-}f_\epsilon 
\label{eq:constr3} \\
&=&\widehat{T}_{i_0}\cdots \widehat{T}_{n-1}\widehat{T}_{w_2}f_\epsilon.
\nonumber
\ena
Combining \eqref{eq:constr1}, \eqref{eq:constr2} and \eqref{eq:constr3}, 
we get 
\begin{eqnarray*}
\widehat{T}_{w_2}f_\epsilon&=&(\pm t^{\pm1/2})^{-d_{\epsilon_{n}}+1}\chi_{i_0}^{-1}\, 
\omega\,\widehat{T}_{w_1}f_{\epsilon_n,\epsilon_1,\ldots,\epsilon_{n-1}}.
\end{eqnarray*}
Since $\omega \widehat{T}_{w_1}=\widehat{T}_{w_2}\omega$,
the third property holds.
\end{proof}

\section{Construction of special solutions}\label{sec:non-symm-Mac}

{}From the result in the foregoing sections 
we can construct special solutions of the qKZ equation as follows. 
Find a solution $E$ to the eigenvalue problem \eqref{eq:equiv1} and \eqref{eq:equiv2}. 
Setting $f_{\epsilon}=(\widehat{T}_{w_{\epsilon}^{-}})^{-1}E$,  
we obtain a qKZ family ${\bf f}=\{f_{\epsilon}\}$ of sign $(\pm)$ 
according to the sign $\pm$ in the right hand side of \eqref{eq:equiv2}. 
Define the parameters $\alpha, \, \beta$ and take a function $h(z)$ 
as explained in Section \ref{subsec:qKZ-family}. 
Using these ingredients above we define $G(z_{1}, \ldots , z_{n})$ 
by the formula \eqref{eq:qKZ-solution}. 
Then from Proposition \ref{prop:qKZ-solution} 
$G$ is a solution of the qKZ equation whose parameter 
$q=\pm t^{\pm1/2}$.
Thus the first step of our construction is to solve 
the eigenvalue problem \eqref{eq:equiv1} and \eqref{eq:equiv2}, 
and we can find a solution in terms of non-symmetric Macdonald polynomials. 

In the following we use the wording 
``the eigenvalue problem of sign $(\pm)$'' to refer 
the eigenvalue problem \eqref{eq:equiv1} and \eqref{eq:equiv2} 
where the sign in the right hand side is $\pm$, 
respectively.   

\subsection{Non-symmetric Macdonald polynomials}

For $\lambda=(\lambda_{1}, \ldots , \lambda_{n}) \in \Z^n$, 
we set $z^\lambda=z_1^{\lambda_1}\cdots z_n^{\lambda_n}$.
We introduce the dominance order $\geq$ on the set $\Z^n$:
\begin{eqnarray*}
\lambda\geq\mu \quad \underset{\tiny \rm def}{\Leftrightarrow} \quad 
{\textstyle \sum_{i=1}^{j}\lambda_i\geq\sum_{i=1}^{j}\mu_i}
\, \mbox{ for any $1\leq j\leq n$},
\end{eqnarray*}
and a partial order $\succ$:
\begin{eqnarray*}
\lambda\succ\mu \quad \underset{\tiny \rm def}{\Leftrightarrow} \quad 
\mbox{ $\lambda^+>\mu^+$ \, or \,\,``$\lambda^+=\mu^+$ and $\lambda>\mu$''}.
\end{eqnarray*}

\begin{defn}
For $\lambda=(\lambda_{1}, \ldots , \lambda_{n}) \in \Z^n$,
the non-symmetric Macdonald polynomial 
$E_\lambda=E_{\lambda}(z_{1}, \ldots , z_{n}; t, p)$
with two parameters $t$ and $p$ is a Laurent polynomial
satisfying
\begin{eqnarray}
\widehat{Y}_iE_\lambda&=&t^{\rho(\lambda)_i}p^{\lambda_i}E_\lambda 
\label{eq:eigenY} \\
E_\lambda&=&z^\lambda+\sum_{\mu\prec\lambda} c_\mu z^\mu
\nonumber 
\end{eqnarray}
where $\rho(\lambda):=w_\lambda^+\rho$,
$\rho:=(\frac{n-1}{2},\frac{n-3}{2},\ldots,-\frac{n-1}{2})$.
\end{defn}

Let us recall the action of $\widehat{T}_i$ on $E_\lambda$ following \cite{Kas}.
Put 
\begin{eqnarray*}
f_i(\lambda):=t^{\rho(\lambda)_{i+1}-\rho(\lambda)_{i}}p^{\lambda_{i+1}-\lambda_{i}}. 
\end{eqnarray*}
If $\lambda_i<\lambda_{i+1}$, then
\begin{eqnarray*}
\widehat{T_i}E_\lambda=t^{1/2}E_{s_i\lambda}-\frac{t^{1/2}-t^{-1/2}}{f_i(\lambda)-1}E_\lambda.
\end{eqnarray*}
If $\lambda_i=\lambda_{i+1}$, then
\begin{eqnarray}
\widehat{T_i}E_\lambda=t^{1/2}E_\lambda.
\label{eq:eigen+}
\end{eqnarray}
If $\lambda_i>\lambda_{i+1}$, then
\begin{eqnarray}
\widehat{T_i}E_\lambda=t^{-1/2}
\frac{(tf_i(\lambda)-1)(t^{-1}f_i(\lambda)-1)}
{(f_i(\lambda)-1)^2}
E_{s_i\lambda}-\frac{t^{1/2}-t^{-1/2}}{f_i(\lambda)-1}E_\lambda.
\label{eq:eigen-}
\end{eqnarray}

The parameters $t$ and $p$ are called {\it generic} if
\begin{eqnarray*}
t^lp^m\neq 1 \quad \hbox{for any} \,\, 0\leq l\leq n \,\, \hbox{and} \,\, 0\leq m. 
\end{eqnarray*}
For generic parameters, $E_\lambda$ is well-defined for any $\lambda\in\Z^n$.

\subsection{Generic case}

First we consider the case where the parameters $t$ and $p$ are generic. 
{}From the properties \eqref{eq:eigenY} and \eqref{eq:eigen+}, 
the non-symmetric polynomials give solutions of 
the eigenvalue problem of sign $(+)$. 
Hence we can get solutions of the qKZ equation: 

\begin{prop}
Suppose that the parameters $t$ and $p$ are generic. 
Let $d_{0}, \ldots , d_{N-1}$ be positive integers satisfying 
$\sum_{j=0}^{N-1}d_{j}=n$ 
and set $\delta=(0^{d_0},1^{d_1},\cdots,(N-1)^{d_{N-1}})$. 
Take $\lambda\in\Z^n$ such that $\lambda_i=\lambda_{i+1}$ if $\delta_i=\delta_{i+1}$.
Then the non-symmetric Macdonald polynomial $E_{\lambda}$ is a solution 
of the eigenvalue problem of sign $(+)$, and 
we obtain a solution of the qKZ equation from it by setting $t^{1/2}=q$. 
The parameters $\kappa_{1}, \ldots , \kappa_{N-1}$ in the qKZ equation 
are determined by \eqref{eq:def-kappa} from the exponents 
\begin{eqnarray*}
c_{i}=q^{d_{i}-1+2\rho(\lambda)_{d_{0}+\cdots +d_{i}}}p^{\lambda_{d_{0}+\cdots +d_{i}}}.   
\end{eqnarray*}
\end{prop}

We note that in the case where $d_{0}=\cdots=d_{N-1}=1$ 
the requirement \eqref{eq:equiv2} becomes empty. 
Hence any non-symmetric Macdonald polynomial is also a solution to 
the eigenvalue problem of sign $(-)$ in this special case, 
and we obtain the following proposition:  

\begin{prop}\label{prop:minus-generic}
Suppose that the parameters $t$ and $p$ are generic. 
In the case where $d_{0}=\cdots=d_{N-1}=1$, and hence $n=N$, 
any non-symmetric Macdonald polynomial $E_{\lambda}$ 
creates a solution of the qKZ equation. 
The parameters $\kappa_{1}, \ldots , \kappa_{N-1}$ are 
determined by \eqref{eq:def-kappa}, 
where $c_i=(-1)^{n-1}q^{-2\rho(\lambda)_{i+1}}p^{\lambda_{i+1}}$. 
\end{prop}

\begin{proof}
Here we only give the calculation of the exponents. 
{}From Theorem \ref{thm:constr} we have 
$c_{i}=t^{\rho(\lambda)_{d_{0}+\cdots +d_{i}}}p^{\lambda_{d_{0}+\cdots +d_{i}}}$. 
To construct the solution of the qKZ equation we set $t^{1/2}=-q^{-1}$. 
Since $2\rho(\lambda)_{i} \equiv n-1 \, ({\rm mod} 2)$ for all $i$, we find
\begin{eqnarray*}
&& 
c_{i}=(-q^{-1})^{2\rho(\lambda)_{d_{0}+\cdots +d_{i}}}p^{\lambda_{d_{0}+\cdots +d_{i}}}=
(-1)^{n-1}q^{-2\rho(\lambda)_{d_{0}+\cdots +d_{i}}}p^{\lambda_{d_{0}+\cdots +d_{i}}} \\ 
&& \quad {}=
(-1)^{n-1}q^{-2\rho(\lambda)_{i+1}}p^{\lambda_{i+1}}. 
\end{eqnarray*} 
\end{proof}

\begin{rem}\label{rem:branch}
When we determine $\kappa_{j}$'s and $\alpha$ 
by \eqref{eq:def-kappa} and \eqref{eq:def-alphabeta} in practice, 
the branch of $(\prod_{i}c_{i})^{1/N}$ should be settled suitably. 
In the situation described in Proposition \ref{prop:minus-generic} 
the exponent $c_{i}$ is a value of the form $(-1)^{n-1}\tilde{c}_{i}$, where 
$\tilde{c}_{i}$ is a positive real number. 
Then we set $(\prod_{i}c_{i})^{1/N}=(-1)^{n-1}(\prod_{i}\tilde{c}_{i})^{1/N}$ 
and determine $\kappa_{j}$'s and $\alpha$. 
In Theorem \ref{thm:specialized} below the situation is the same, 
and we take the same branch. 
\end{rem}

\subsection{Specialized case}\label{subsec:specialized}

In Proposition \ref{prop:minus-generic} 
we saw that 
any non-symmetric Macdonald polynomial gives a solution to 
the eigenvalue problem of sign $(-)$, 
but it is in the very special case. 
In order to solve this problem in general, 
we need to find an eigenfunction of the Demazure-Lusztig operator 
with the eigenvalue $-t^{-1/2}$, 
and this is not the situation in \eqref{eq:eigen+}. 
However, if $f_{i}(\lambda)=t$ in \eqref{eq:eigen-}, 
then $E_{\lambda}$ becomes such an eigenfunction. 
It should be noted that the relation $f_{i}(\lambda)=t$ implies that 
the parameters $t$ and $p$ are not generic. 
In the rest of this paper we consider this kind of case. 

Let $k$ and $r$ be integers such that 
$1\leq k\leq {\rm min}\{n-1, N\}, \, r \ge 2$, and $k+1$ and $r-1$ are coprime.   
We assume that $t,p$ are not roots of unity
and take a specialization $t^{k+1}p^{r-1}=1$.
To be more precise we specialize $t$ and $p$ as follows:
\begin{eqnarray}
t=u^{r-1}, && p=u^{-(k+1)}
\label{eq:spec}
\end{eqnarray}
where $u$ is not root of unity. 
We will set $q=-t^{-1/2}$ and take $u=q^{-\frac{2}{r-1}}$. 
Then we have $p=q^{\frac{2(k+1)}{r-1}}$ and 
the level of the qKZ equation is equal to $\frac{k+1}{r-1}-N$.  

We call $\lambda=(\lambda_{1}, \ldots , \lambda_{n}) \in \Z^n$ {\it admissible} if
\begin{eqnarray*}
&& 
\lambda^+_i-\lambda^+_{i+k}\geq r-1 \quad \hbox{for any $1\leq i\leq n-k$, and} \\ 
&& 
\lambda^+_i-\lambda^+_{i+k}=r-1 \quad \hbox{only if} \,\, w_\lambda^+(i)<w_\lambda^+(i+k). 
\end{eqnarray*}
The following statement is a corollary of Theorem 3.11 in \cite{Kas}:
\begin{lem}\label{lem:specialized}
For any admissible $\lambda\in\Z^n$, 
the non-symmetric Macdonald polynomial $E_\lambda$ is well-defined
under the specialization $(\ref{eq:spec})$.
If $\lambda\in\Z^n$ is admissible and $s_i\lambda$
is not admissible, then $\widehat{T_i}E_\lambda=-t^{-1/2}E_\lambda$.
\end{lem}

Let $m$ and $l$ be integers satisfying
$n=km+l$ and $0\leq l\leq k-1$.
Let $(d^{(0)},\ldots,d^{(k-1)})$ be a permutation of $((m+1)^l,m^{k-l})$.
Note that $\sum_{j=0}^{k-1}d^{(j)}=n$. 
Take a dominant element $a=(a_{1}, \ldots , a_{k}) \in \Z^k$ satisfying
\begin{eqnarray}
a_1-a_k\leq r-1 \quad \hbox{and} \quad w_{a}^-((m+1)^l,m^{k-l})=(d^{(0)}, \ldots , d^{(k-1)}),  
\label{cond:dominant-a}
\end{eqnarray}
where $w_{a}^{-} \in S_{k}$ (see Definition \ref{def:shortest-element}). 
Now define $\lambda \in \mathbb{Z}^{n}$ by 
\begin{eqnarray*}
&& 
\lambda_i=a_i \quad \hbox{for} \quad 1\leq i \leq k, \quad \hbox{and} \\ 
&& 
\lambda_i-\lambda_{i+k}=r-1 \quad \hbox{for} \quad 1\leq i \leq n-k. 
\end{eqnarray*}
Then $\lambda$ is admissible. 
For simplicity, we write $w=w_a^-$ and define $\mu \in S_n\lambda$ by
\begin{eqnarray*}
\mu&=&(\lambda_{w^{-1}(1)},\lambda_{w^{-1}(1)+k},
\lambda_{w^{-1}(1)+2k},\ldots,\lambda_{w^{-1}(1)+b_{1}k},\\
&&\,\,
\lambda_{w^{-1}(2)},\lambda_{w^{-1}(2)+k},
\lambda_{w^{-1}(2)+2k},\ldots,\lambda_{w^{-1}(2)+b_{2}k},\\
&&\,\, 
\ldots,\\
&&\,\, 
\lambda_{w^{-1}(k)},\lambda_{w^{-1}(k)+k},
\lambda_{w^{-1}(k)+2k},\ldots,\lambda_{w^{-1}(k)+b_{k}k}), 
\end{eqnarray*}
where $b_{j}:=m-\theta(w^{-1}(j)>l)$ (see \eqref{eq:def-theta} for the definition of $\theta(P)$).
\bigskip 

\noindent{\it Example.}
Set $n=13$ and $k=5$, 
and consider the case of $(d^{(0)}, d^{(1)}, d^{(2)}, d^{(3)}, d^{(4)})=(3,2,2,3,3)$.  
Then the condition \eqref{cond:dominant-a} 
for a dominant $a=(a_{1}, \ldots , a_{5}) \in \mathbb{Z}^{5}$ implies that 
$a_{1} \ge a_{2} > a_{3}=a_{4}=a_{5}$. 
Now suppose that $r=6$ and take $a=(13,10,9,9,9)$. 
Then $\lambda$ and $\mu$ are given by 
\begin{eqnarray*} 
&& 
\lambda=(13,10,9,9,9,8,5,4,4,4,3,0,-1), \\
&& 
\mu=(9,4,-1,9,4,9,4,10,5,0,13,8,3). 
\end{eqnarray*}
\medskip 

Now let $(d_0,\ldots,d_{N-1})$ be a subdivision of $(d^{(0)},\ldots,d^{(k-1)})$, 
that is, $d_{i}>0$ and 
$d_{i_j}+\cdots+d_{i_{j+1}-1}=d^{(j)}$ for some $0=i_{0}<i_1<\cdots <i_{k-1}<i_{k}=N$.
It is easy to see that $\mu$ is also admissible
and $s_i\mu$ is not admissible if $\delta_i=\delta_{i+1}$, 
where $\delta=(0^{d_{0}},1^{d_{1}}, \cdots, (N-1)^{d_{N-1}})$.
{}From Lemma \ref{lem:specialized}, 
$E_\mu$ is a solution of the eigenvalue problem of sign $(-)$.  
Therefore we get the following theorem
(see Remark \ref{rem:branch}). 

\begin{thm}\label{thm:specialized}
The non-symmetric Macdonald polynomial $E_{\mu}$ 
with the specialization \eqref{eq:spec} and $t^{1/2}=-q^{-1}$ 
creates a solution of the qKZ equation of level $\frac{k+1}{r-1}-N$. 
The parameters $\kappa_{1}, \ldots , \kappa_{N-1}$ are determined 
by \eqref{eq:def-kappa} from the exponents 
$c_{i}=(-1)^{n-1}q^{A_{i}}$, where 
\begin{eqnarray*}
A_{i}:=d_i-1-2\rho(\mu)_{d_{0}+\cdots+d_i}+\frac{2(k+1)}{r-1}\mu_{d_{0}+\cdots+d_i}. 
\end{eqnarray*}
\end{thm}

\section{The matrix element of the vertex operators}\label{sec:VO}

Here we see that our special solutions constructed in Theorem \ref{thm:specialized} 
contain the matrix element of the vertex operators 
in the case where $k=N$ and $r=2$.  

First we recall the definition and the properties 
of the vertex operators following \cite{DO}. 
Let $\Lambda_{l} \, (l=0, \ldots , N-1)$ be the fundamental weights of 
$U_{q}(\widehat{sl}_{N})$.
The symmetric bilinear form $( \cdot | \cdot )$ on the weight lattice is defined by 
\begin{eqnarray*}
(\Lambda_{i}|\Lambda_{j})=\frac{i(N-j)}{N} \qquad (i\le j).   
\end{eqnarray*}
Hereafter the index $l$ in $\Lambda_{l}$ should read modulo $N$. 
The weight of $v_{\epsilon} \in V$ is given by 
${\rm wt}v_{\epsilon}=\bar{\Lambda}_{\epsilon+1}-\bar{\Lambda}_{\epsilon}$, 
where $\bar{\Lambda}_{l}:=\Lambda_{l}-\Lambda_{0}$. 

Denote by $V(\Lambda_{l})$ the irreducible integrable highest weight module 
of $U_{q}(\widehat{sl}_{N})$ 
with the highest weight $\Lambda_{l}$, 
and by $|\Lambda_{l} \rangle $ its highest weight vector. 
Let $\langle \Lambda_{l} | \in V(\Lambda_{l})^{*}$ be the dual vector 
satisfying $\langle \Lambda_{l} | \Lambda_{l} \rangle=1$.  

The vertex operator $\widetilde{\Phi}^{(l)}(z) \, (0 \le l \le N-1)$ of type I 
\footnote{This operator $\widetilde{\Phi}^{(l)}(z)$ is equal to 
$\widetilde{\Phi}_{\Lambda_{l+1}}^{\Lambda_{l}, V^{(1)}}(z)$ in \cite{DO}.}  
is the intertwiner 
\begin{eqnarray*}
\widetilde{\Phi}^{(l)}(z): \, 
V(\Lambda_{l+1}) \rightarrow V(\Lambda_{l}) \otimes V 
\end{eqnarray*}
normalized as 
$\widetilde{\Phi}^{(l)}(z) | \Lambda_{l+1} \rangle=|\Lambda_{l} \rangle\otimes v_{l}+ \cdots$.  
In the following we often omit the upper index $(l)$ of the vertex operators. 
Write the operator $\widetilde{\Phi}(z)$ as 
$\widetilde{\Phi}(z)(\cdot)=\sum_{\epsilon=0}^{N-1}\widetilde{\Phi}_{\epsilon}(z)(\cdot) \otimes v_{\epsilon}$.  
Then the following commutation relation holds:
\begin{eqnarray}
\widetilde{\Phi}_{\epsilon_{2}}(z_{2})\widetilde{\Phi}_{\epsilon_{1}}(z_{1})=
\kappa(z)\sum_{\epsilon_{1}', \epsilon_{2}'=1}^{N}
\bar{R}(z_{1}/z_{2})_{\epsilon_{1}, \epsilon_{2}}^{\epsilon_{1}', \epsilon_{2}'}
\widetilde{\Phi}_{\epsilon_{1}'}(z_{1})\widetilde{\Phi}_{\epsilon_{2}'}(z_{2}), 
\label{eq:commrel-VO}
\end{eqnarray}
where 
\begin{eqnarray*}
\kappa(z)=z^{\frac{1}{N}-1}
\frac{(q^{2N}z^{-1}; q^{2N})_{\infty} (q^{2}z; q^{2N})_{\infty}}
{(q^{2N}z; q^{2N})_{\infty} (q^{2}z^{-1}; q^{2N})_{\infty}}.   
\end{eqnarray*}

Set $\Phi^{(l)}(z)=z^{\Delta_{l}-\Delta_{l+1}}\widetilde{\Phi}^{(l)}(z)$, 
where $\Delta_{l}:=\frac{l(N-l)}{2N}$.
Consider the matrix element 
\begin{eqnarray*}
G^{(ij)}(z_{1}, \ldots , z_{n})=
\langle \Lambda_{i} | \Phi(z_{1}) \cdots \Phi(z_{n}) | \Lambda_{j} \rangle \in V^{\otimes n}.
\end{eqnarray*}

\begin{thm}\cite{FR}\label{thm:qKZ}
The function $G^{(ij)}$ satisfies the qKZ equation of level one. 
The parameters $\kappa_{1}, \ldots , \kappa_{N-1}$ are determined by 
$\prod_{l=1}^{N-1}\kappa_{l}^{h_{l}}=q^{-\bar{\Lambda}_{i}-\bar{\Lambda}_{j}-2\bar{\rho}}$,  
where $\bar{\rho}:=\sum_{l=1}^{N-1}\bar{\Lambda}_{l}$. 
\end{thm}

{}From the definition of the vertex operators 
we have $G^{(ij)}=0$ unless $i-j+n \equiv 0 \, ({\rm mod}\, N)$.
In the following we assume that $n\ge N$ and 
$i-j+n \equiv 0 \, ({\rm mod}\, N)$. 
The function $G^{(ij)}$ is expanded as 
\begin{eqnarray*}
&& 
G^{(ij)}(z_{1}, \ldots , z_{n})=
\sum_{\epsilon_{1}, \ldots , \epsilon_{n}} 
G^{(ij)}_{\epsilon_{1}, \ldots , \epsilon_{n}}(z_{1}, \ldots , z_{n})\, 
v_{\epsilon_{1}} \otimes \cdots \otimes v_{\epsilon_{n}}, \\ 
&& 
G^{(ij)}_{\epsilon_{1}, \ldots , \epsilon_{n}}(z_{1}, \ldots , z_{n}):=
\langle \Lambda_{i} | \Phi_{\epsilon_{1}}(z_{1}) \cdots 
\Phi_{\epsilon_{n}}(z_{n}) | \Lambda_{j} \rangle.
\end{eqnarray*} 
Then $G_{\epsilon}^{(ij)}=G_{\epsilon_{1}, \ldots , \epsilon_{n}}^{(ij)}=0$ unless 
$\epsilon \in I_{d_{0}, \ldots , d_{N-1}}$, 
where $d_{0}, \ldots , d_{N-1}$ are determined by $\sum_{l=0}^{N-1}d_{l}=n$ and 
\begin{eqnarray}
d_{l}=d_{l-1}+\delta_{l,i}-\delta_{l,j} \qquad (0 \le l \le N-1).  
\label{eq:def-d}
\end{eqnarray}
Here the index $l$ should read modulo $N$. 

The following formula is due to Nakayashiki \cite{N}: 
\begin{prop}
Set $\delta=(0^{d_{0}}, \ldots , (N-1)^{d_{N-1}})$. Then we have 
\begin{eqnarray*}
G^{(ij)}_{\delta}(z_{1}, \ldots , z_{n})=
K^{(j)}(z_{1}, \ldots , z_{n})
\prod_{a=1}^{d_{0}+\cdots +d_{j-1}}z_{a}^{-1} 
\prod_{1 \le a<b \le n \atop \delta_{a}=\delta_{b}}(z_{a}-q^{2}z_{b}).   
\end{eqnarray*} 
Here $K^{(j)}$ is defined by 
\begin{eqnarray}
&& 
K^{(j)}(z_{1}, \ldots , z_{n}) 
\label{eq:def-K-VO} \\ 
&& {}:=
c_{ij}^{(n)}\, 
\prod_{a=1}^{n}z_{a}^{\frac{1}{N}a+\frac{1}{2N}(-2n+2j+N-1)} \, 
\prod_{1 \le a<b \le n}
\frac{(q^{2N+2}z_{b}/z_{a}; q^{2N})_{\infty}}{(q^{2N}z_{b}/z_{a}; q^{2N})_{\infty}},  
\nonumber
\end{eqnarray}
where $c_{ij}^{(n)}$ is a certain constant.
\end{prop}

We have the following proposition. 

\begin{prop}
Set $F^{(ij)}_{\epsilon}:=G_{\epsilon}^{(ij)}/K^{(j)}$ for $\epsilon \in I_{d_{0}, \ldots , d_{N-1}}$. 
Then $\{F^{(ij)}_{\epsilon}\}$ is a qKZ family of sign $(-)$ whose parameter $t^{1/2}=-q^{-1}$. 
The exponents are given by 
\begin{eqnarray}
c_{\epsilon}=(-1)^{n-1}
q^{n-2N-2j+2\epsilon-\theta(\epsilon<i)-\theta(\epsilon<j)+\frac{i-j+n}{N}}. 
\label{eq:exponents-VO}
\end{eqnarray}   
\end{prop}

\begin{proof}
Set $F^{(ij)}=G^{(ij)}/K^{(j)}$. 
{}From the commutation relation \eqref{eq:commrel-VO} and 
the definition \eqref{eq:def-K-VO} of $K^{(j)}$, 
we have 
\begin{eqnarray}
&& 
F^{(ij)}( \ldots , z_{a+1}, z_{a}, \ldots ) 
\label{eq:commrel-F} \\ 
&& {}=
\frac{z_{a+1}-q^{2}z_{a}}{z_{a}-q^{2}z_{a+1}} P_{a, a+1} \, 
\bar{R}_{a,a+1}(z_{a}/z_{a+1})
F^{(ij)}(\ldots , z_{a}, z_{a+1}, \ldots), 
\nonumber 
\end{eqnarray} 
where $P$ is the transposition $P(u \otimes v):=v \otimes u$. 
This is equivalent to the first and second properties
in Definition \ref{def:qKZ-family}. 

{}From \eqref{eq:commrel-F} and the qKZ equation for $G^{(ij)}$ 
we see that 
\begin{eqnarray*}
F^{(ij)}(q^{2(N+1)}z_{1}, z_{2}, \ldots , z_{n})&=& 
(-1)^{n-1}q^{n-N-2j-1+(n-2j)/N} \, 
(q^{-\bar{\Lambda}_{i}-\bar{\Lambda}_{j}-2\bar{\rho}})_{1} \\ 
&\times& 
P_{1,2} \cdots P_{n-1,n}F^{(ij)}(z_{2}, \ldots , z_{n}, z_{1}). 
\end{eqnarray*}
It is easy to derive the third property in Definition \ref{def:qKZ-family} from this equality. 
\end{proof}

For the exponents \eqref{eq:exponents-VO} and $p=q^{2(N+1)}$, 
the parameters $\alpha$ and $\beta$ are determined by \eqref{eq:def-alphabeta} as follows: 
\begin{eqnarray*}
\alpha=\frac{1}{2N}(-2n+2j+N-1), \quad 
\beta=\frac{1}{N}.  
\end{eqnarray*}
Moreover the function $h(z):=(q^{2N+2}z; q^{2N})_{\infty}/(q^{2N}z; q^{2N})_{\infty}$ is 
a solution to the difference equation \eqref{eq:def-h}. 
Thus the function $K^{(j)}$ is restored from the ingredients above 
by the definition \eqref{eq:normalization-factor} 
up to constant multiplication. 
 
Let us check that the matrix element $G^{(ij)}$ is contained in 
our special solutions constructed in Theorem \ref{thm:specialized} 
in the case where $k=N$ and $r=2$. 
Since the function $K^{(j)}$ is restored as described above, 
it suffices to show that the extremal component 
\begin{eqnarray*}
F_{\delta}^{(ij)}=
\prod_{a=1}^{d_{0}+\cdots +d_{j-1}}z_{a}^{-1} 
\prod_{1 \le a<b \le n \atop \delta_{a}=\delta_{b}}(z_{a}-q^{2}z_{b})  
\end{eqnarray*}
is equal to the non-symmetric Macdonald polynomial $E_{\mu}$
for a suitable $\mu \in \mathbb{Z}^{n}$ with the specialization $t^{N+1}p=1$.  

Recall that $d_{l} \, (0 \le l \le N-1)$ is defined by \eqref{eq:def-d}. 
More explicitly we have 
\begin{eqnarray*}
d:=(d_{0}, \ldots , d_{N-1})=
\left\{ 
\begin{array}{ll}
(m^{i}, (m+1)^{j-i}, m^{N-j}) & \hbox{if $i\le j$},  \\
((m+1)^{j}, m^{i-j}, (m+1)^{N-i}) & \hbox{if $i > j$}. 
\end{array}
\right.
\end{eqnarray*}
Here $m:=[n/N]$.  
Now set $a=(a_{1}, \ldots , a_{N}) \in \mathbb{Z}^{N}$ by 
\begin{eqnarray*}
a=\left\{ 
\begin{array}{ll}
((m-1)^{N-i}, (m-2)^{i}) & \hbox{if $i\le j$},  \\
(m^{N-i}, (m-1)^{i}) & \hbox{if $i > j$}. 
\end{array}
\right. 
\end{eqnarray*}
It satisfies $w_{a}^{-}(d^{+})=d$. 
This element $a$ determines $\mu \in \mathbb{Z}^{n}$ as described in Section \ref{subsec:specialized}, 
and it is given by 
\begin{eqnarray*}
\mu=\left\{ 
\begin{array}{ll}
((m-2, \ldots , -1)^{i}, (m-1, \ldots , -1)^{j-i}, (m-1, \ldots , 0)^{N-j}) & \hbox{if $i\le j$},  \\
((m-1, \ldots , -1)^{j}, (m-1, \ldots , 0)^{i-j}, (m, \ldots , 0)^{N-i}) & \hbox{if $i > j$}. 
\end{array}
\right.  
\end{eqnarray*} 

\begin{prop}\label{prop:difference-product}
For $\mu \in \mathbb{Z}^{n}$ defined above, 
$E_{\mu}$ is equal to $F_{\delta}^{(ij)}$ with the specialization $t^{N+1}p=1$ and $t^{1/2}=-q^{-1}$.  
Therefore the solution of the qKZ equation determined from $E_{\mu}$ 
coincides with the matrix element $G^{(ij)}$ of 
the vertex operators up to constant multiplication. 
\end{prop}

To prove the proposition we use the following lemma obtained in \cite{Kas}: 

\begin{lem}\label{lem:wheel}
Suppose $1\leq N \leq n-1$.
For a Laurent polynomial $f$ we call
the following vanishing property the wheel condition:
$f(z_1, \ldots, z_n)=0$ if $z_{i_1}=t^{-1}z_{i_2}=\cdots=t^{-N}z_{i_{N+1}}$
for any $1\leq i_1<\cdots<i_{N+1}\leq n$.
Then the set of the non-symmetric Macdonald polynomials 
\begin{eqnarray*}
\{ E_{\lambda} \, | \, \hbox{$\lambda$ is admissible} \} 
\end{eqnarray*}
with the specialization $t^{N+1}p=1$ forms
 a basis of the space of all Laurent polynomials satisfying the wheel condition.
\end{lem}

\begin{proof}[Proof of Proposition \ref{prop:difference-product}]
In the case of $n=N$ the equality $F_{\delta}^{(ij)}=E_{\mu}$ follows from 
the definition of $E_{\mu}$. 
Let us consider the case of $n>N$. 
Note that $F_{\delta}^{(ij)}$ satisfies the wheel condition and 
the top term in $F_{\delta}^{(ij)}$ is equal to $x^\mu$.  
Hence, from Lemma \ref{lem:wheel},  
$F_{\delta}^{(ij)}$ is a linear combination of
$E_\nu$ where $\nu$ is admissible, $\nu\preceq\mu$ and $|\nu|=|\mu|$.
Let $\nu$ be such an element.
By definition, $\nu^+<\mu^+$ or ``$\nu^+=\mu^+$ and $\nu\leq\mu$".
If $\nu^+<\mu^+$, then $\nu^+_1=\mu^+_1,\ldots,\nu^+_{s-1}=\mu^+_{s-1}$,
and $\nu^+_s<\mu^+_s$ for some $s$.
Since $\nu$ is admissible, $\nu^+$ is also admissible.
However, the admissibility of $\nu^+$ and the definition of $\mu$ 
impose an upper bound for the other components $\nu^+_{s'}\leq\mu^+_{s'}$ ($s'>s$).
It is inconsistent with the condition $|\nu^+|=|\mu^+|$.
Thus $\nu$ must satisfy $\nu^+=\mu^+$ and $\nu\leq\mu$.
If $\nu<\mu$, then
$\nu_1=\mu_1,\ldots,\nu_{s-1}=\mu_{s-1}$, and $\nu_s<\mu_s$ for some $s$.
However it is not compatible with the admissibility of $\nu$.
Hence the only possibility is $\nu=\mu$.
\end{proof}

\bigskip 
\noindent{\it Acknowledgements.} 
Research of YT is supported by Grant-in-Aid for 
Young Scientists (B) No.\,17740089. 
Research of MK is supported by Grant-in-Aid for 
JSPS Fellows No. 17-2106. 
We are deeply grateful to Saburo Kakei and Yoshihisa Saito for valuable discussions. 
MK also thanks Vincent Pasquier for interests. 
YT thanks Takeshi Ikeda and Hiroshi Naruse for discussions.

\end{document}